\documentclass[12pt]{amsart}

\usepackage{fullpage,url}
\usepackage{amsmath}
\usepackage{amsfonts}
\usepackage{amssymb}
\usepackage{amscd}

\newcommand{\pp}{{\mathfrak p}}

\DeclareMathOperator{\h}{h}

\DeclareMathOperator{\tor}{tor}
\DeclareMathOperator{\alg}{alg}
\DeclareMathOperator{\hhat}{\widehat{h}}

\DeclareMathOperator{\Frac}{Frac}

\newtheorem{theorem}{Theorem}[section]
\newtheorem{lemma}[theorem]{Lemma}

\newtheorem{proposition}[theorem]{Proposition}

\theoremstyle{definition}
\newtheorem{definition}[theorem]{Definition}
\newtheorem{Claim}[theorem]{Claim}

\newtheorem{conjecture}[theorem]{Conjecture}
\newtheorem{example}[theorem]{Example}

\theoremstyle{remark}
\newtheorem{remark}[theorem]{Remark}

\title{The Tate-Voloch conjecture for Drinfeld modules}
\author{Dragos Ghioca}

\begin{document}
\begin{abstract}
We study the $v$-adic distance from the torsion of a Drinfeld module to an affine variety. 
\end{abstract}
\maketitle
\section{Introduction}
\label{se:intro}
\footnotetext[1]{2000 AMS Subject Classification: Primary, 11G09; Secondary, 11G50}

For a semi-abelian variety $S$ and an algebraic subvariety $X\subset S$, the Manin-Mumford conjecture characterizes the subset of torsion points of $S$ contained in $X$. The Tate-Voloch conjecture characterizes the distance from $X$ of a torsion point of $S$ not contained in $X$.

Let $\mathbb{C}_p$ be the completion of a fixed algebraic closure $\mathbb{Q}_p^{\alg}$ of $\mathbb{Q}_p$. Let $\lambda(\cdot,X)$ be the $p$-adic proximity to $X$ function as defined in \cite{Sca1} (see also our definition of $v$-adic distance to an affine subvariety). Tate and Voloch conjectured:

\begin{conjecture}[Tate,Voloch]
\label{C:cvt}
Let $G$ be a semi-abelian variety over $\mathbb{C}_p$. Let $X\subset G$ be a subvariety defined over $\mathbb{C}_p$. Then there is a constant $N\in\mathbb{N}$ such that for any torsion point $\zeta\in G(\mathbb{C}_p)$ either $\zeta\in X$ or $\lambda_p(\zeta,X)\le N$.
\end{conjecture}

The above conjecture was proved by Thomas Scanlon for all semi-abelian varieties defined over $\mathbb{Q}_p^{\alg}$ (see \cite{Sca1} and \cite{Sca2}).

In this paper we prove two Tate-Voloch type theorems for Drinfeld modules. Our motivation is to show that yet another question for semi-abelian varieties has a counterpart for Drinfeld modules (see \cite{Sca3} for a version of the Manin-Mumford theorem for Drinfeld modules of generic characteristic and see \cite{IMRN} for a version of the Mordell-Lang theorem for Drinfeld modules).

In Section~\ref{S:all} we state our results. Our first result (Theorem~\ref{T:tv}) shows that if a torsion point of a Drinfeld module $\phi:A\rightarrow K\{\tau\}$ is close $w$-adically to a variety $X$ with respect to all places $w$ extending a fixed place $v$ of the ground field $K$, then the torsion point lies on $X$. We prove Theorem~\ref{T:tv} in Section~\ref{S:one}. Our bound for how "close $w$-adically to $X$" means "lying on $X$" is effective. Our second result (Theorem~\ref{T:Mattuck}) refers to proximity with respect to one fixed extension of a place $v$ of $K$. We will prove Theorem~\ref{T:Mattuck} in Section~\ref{S:future}. We also note that due to the fact that in Theorem~\ref{T:Mattuck} we work with a fixed extension of a place of $K$, there is a different normalization for the valuation we are working as opposed to the setting in Theorem~\ref{T:tv}.

\section{Statement of our main results}
\label{S:all}

Before stating our results we introduce the definition of a Drinfeld module (for more details, see \cite{teza}).

Let $p$ be a 
prime number and let $q$ be a power of $p$.
We let $C$ be a nonsingular projective curve defined
over $\mathbb{F}_q$ and we 
fix a closed point $\infty$ on $C$. Then we define $A$ as the ring of functions
on $C$ that are regular 
everywhere except possibly at $\infty$.

We let $K$ be a field extension of $\mathbb{F}_q$ and we fix an algebraic closure of $K$, denoted
$K^{\alg}$. We fix a 
morphism $i:A\rightarrow K$. We define the operator
$\tau$ as the power of the 
usual Frobenius with the property that for every $x\in
K^{\alg}$, $\tau(x)=x^q$. 
Then we let $K\{\tau\}$ be the ring of polynomials in
$\tau$ with coefficients 
in $K$ (the addition is the usual one, while the multiplication is the composition of functions).

A Drinfeld module over $K$ is a ring morphism $\phi:A\rightarrow
K\{\tau\}$ for which the 
coefficient of $\tau^0$ in $\phi_a$ is $i(a)$ for
every $a\in A$, and there 
exists $a\in A$ 
such that $\phi_a\ne i(a)\tau^0$. We call $\phi$ a Drinfeld module of
generic characteristic 
if $\ker(i)=\{0\}$ and we call $\phi$ a Drinfeld
module of finite 
characteristic if $\ker(i)\ne \{0\}$. In the generic characteristic case we assume $i$ extends to an embedding of $\Frac(A)$ (which is the function field of the projective nonsingular curve $C$) into $K$.

For every nonzero $a\in A$, let the \emph{$a$-torsion} $\phi[a]$ of $\phi$ be the set of all $x\in K^{\alg}$ such that $\phi_a(x)=0$. Let the torsion submodule of $\phi$ be $\bigcup_{a\in A\setminus\{0\}}\phi[a]$. 

For every $g\ge 1$, let $\phi$ act diagonally on $\mathbb{G}_a^g$. An element $(x_1,\dots,x_g)\in (K^{\alg})^g$ is called a torsion element of $\phi$, if for every $i\in\{1,\dots,g\}$, $x_i\in\phi_{\tor}$. 

For each field extension $L$ of $K$ and for each valuation $w$ on $L$ we define the $w$-adic distance to an affine subvariety $X\subset\mathbb{G}_a^g$ defined over $L$. Let $I_X$ be the vanishing ideal in $L[X_1,\dots,X_g]$ of $X$. Let $R_w\subset L$ be the valuation ring of $w$. If $P\in\mathbb{G}_a^g(L)$, then the $w$-adic distance from $P$ to $X$ is
\begin{equation}
\label{E:distance}
\lambda_w(P,X):=\min\{w(f(P))\mid f\in I_X\cap R_w[X_1,\dots,X_g]\}.
\end{equation}

We denote by $M_K$ the set of all discrete valuations on $K$. Similarly, for each field extension $L$ of $K$ we also denote by $M_L$ the set of all discrete valuations on $L$. Finally, we note that unless otherwise stated, each valuation is normalized so that its range is precisely $\mathbb{Z}\cup\{+\infty\}$ (our convention is that the valuation of $0$ is $+\infty$). Our Theorem~\ref{T:tv} is valid for all fields $K$ equipped with a \emph{coherent good set of valuations}. 

\begin{definition}
\label{D:good sets}
We call a subset $U\subset M_K$ equipped with a function 
$d:U\rightarrow\mathbb{R}_{>0}$ a \emph{good set of valuations}
if the following properties are satisfied

(i) for every nonzero $x\in K$, there are finitely
many $v\in U$ such that $v(x)\ne 0$.

(ii) for every nonzero $x\in K$,
$$\sum_{v\in U}d(v)\cdot v(x)=0.$$
The positive real number $d(v)$ will be called the \emph{degree} of the 
valuation $v$. When we say that the positive real number $d(v)$ is associated to 
the valuation $v$, we understand that the degree of $v$ is $d(v)$. 

When $U$ is a good set of valuations, we will refer to
property (ii) as the sum formula for $U$.
\end{definition}

\begin{definition}
\label{D:loccoh}
Let $v\in M_K$ of degree $d(v)$. We say that the valuation $v$ is \emph{coherent} if for every finite extension $L$ of $K$,
\begin{equation}
\label{E:(iii)}
\sum_{\substack{w\in M_{L}\\
w|v}}e(w|v)f(w|v)=[L:K],
\end{equation}
where $e(w|v)$ is the ramification index and $f(w|v)$ is the relative degree 
between the residue field of $w$ and the residue field of $v$.

Condition \eqref{E:(iii)} says that $v$ is \emph{defectless} in $L$. In this 
case, we also let the degree of any $w\in M_L$, $w|v$ be
\begin{equation}
\label{E:(ii)}
d(w)=\frac{f(w|v)d(v)}{[L:K]}.
\end{equation}
\end{definition}

\begin{definition}
\label{D:coherent good sets}
We let $U_K$ be a good set of valuations on $K$. We call $U_K$ a \emph{coherent} 
good set of
valuations if for every $v\in U_K$, the valuation $v$ is coherent.
\end{definition}

\begin{remark}
\label{R:important remark}
Using the argument from page $9$ of \cite{Ser}, we conclude that in Definition~\ref{D:coherent good sets}, if for each finite extension $L$ of $K$ we let $U_L\subset M_L$ be the set of valuations lying above valuations in $U_K$, then $U_L$ is a good set of valuations.
\end{remark} 

\begin{example}
Let $V$ be a projective, regular in codimension $1$ variety defined over a finite field. Then the function field $F$ of $V$ is equipped with a coherent good set of valuations associated to each irreducible divisor of $V$. Hence every finitely generated field is equipped with at least one coherent good set of valuations (different sets of valuations correspond to different projective, regular in codimension $1$ varieties with the same function field). For more details see \cite{Ser} or Chapter $4$ of \cite{teza}.
\end{example}

We prove the following Tate-Voloch type theorem for Drinfeld modules.

\begin{theorem}
\label{T:tv}
Assume $U_K$ is a coherent good set of valuations on $K$ and let $v\in U_K$ have degree $d(v)$. Let $\phi:A\rightarrow K\{\tau\}$ be a Drinfeld module. Let $X\subset\mathbb{G}_a^g$ be a closed $K$-subvariety of the $g$-dimensional affine space.

There exists a constant $C>0$ (depending on $X$ and $d(v)$) such that for every finite extension $L$ of $K$ and for every torsion point $P\in\mathbb{G}_a^g(L)$ of $\phi$, either $P\in X(L)$ or there exists $w\in M_L$ lying over $v$ such that $\lambda_w(P,X)\le C\cdot e(w|v)$.
\end{theorem}

\begin{remark}
\label{R:connection}

There are two significant differences between our Tate-Voloch type theorem and Conjecture~\ref{C:cvt}. We show that a torsion point of the Drinfeld module is on $X$ if it is close to $X$ with respect to all extensions of a fixed valuation $v$ of $K$, not only with respect to one fixed extension of $v$. We will show in Example~\ref{E:infinity} that we cannot always expect proximity of $P$ to $X$ with respect to one fixed extension of $v$ imply that $P$ lies on $X$. The second difference between our Theorem~\ref{T:tv} and Conjecture~\ref{C:cvt} is purely technical. Because we normalized all valuations so that their ranges equal $\mathbb{Z}$, we need to multiply by the corresponding ramification index the constant $C$ in Theorem~\ref{T:tv}.

\end{remark}

\begin{example}
\label{E:infinity}

Let $\phi$ be any Drinfeld module of generic characteristic and let $v_{\infty}$ be a valuation on $K$ extending the valuation on $\Frac(A)$ associated to the closed point $\infty\in C$. We let $K_{\infty}$ be a completion of $K$ with respect to $v_{\infty}$. Then $\phi_{\tor}\subset K_{\infty}^{\alg}$ is not discrete with respect to $v_{\infty}$ (see Section $4.13$ of \cite{Gos}). Hence there exist nonzero torsion points of $\phi$ arbitrarily close to $X:=\{0\}$ in the $v_{\infty}$-adic topology.

\end{example}

For the remaining of Section~\ref{S:all} we fix a valuation $v$ on $K$ (we do not require anymore that $v$ belongs to a good set of valuations on $K$ nor that $v$ is coherent). We let $K_v$ be the completion of $K$ at $v$. We fix an algebraic closure $K_v^{\alg}$ of $K_v$ and extend $v$ to a valuation of $K_v^{\alg}$. In this case, the value group of $v$ is $\mathbb{Q}$. We define as in \eqref{E:distance} the $v$-adic distance from a point $P\in\mathbb{G}_a^g(K_v^{\alg})$ to a fixed affine variety $X$ defined over $K_v^{\alg}$.

Our Theorem~\ref{T:Mattuck} characterizes the distance from $\phi_{\tor}^g$ to a fixed point of $\mathbb{G}_a^g(K_v^{\alg})$. Our theorem is an analogue for Drinfeld modules of a Theorem of Mattuck (see \cite{M}).

\begin{theorem}
\label{T:Mattuck}
Let $\phi:A\rightarrow K\{\tau\}$ be a Drinfeld module. Let $v$ be a place of $K$. If $\phi$ is a Drinfeld module of generic characteristic, then assume $v$ does not lie over the valuation $v_{\infty}$ of $\Frac(A)$, which is associated to the closed point $\infty\in C$. Let $g\ge 1$.

Then for every $Q\in\mathbb{G}_a^g(K_v^{\alg})$ there exists a positive constant $C$ depending on $\phi$, $v$ and $Q$ such that for each $P\in\phi_{\tor}^g$ either $P=Q$ or $\lambda_v(P,Q)<C$.
\end{theorem}

Note that as shown in Example~\ref{E:infinity}, Theorem~\ref{T:Mattuck} does not hold if $v$ extends the place $v_{\infty}$ of $\Frac(A)$, in case $\phi$ has generic characteristic. If $\phi$ has finite characteristic, there is no restriction on $v$ in Theorem~\ref{T:Mattuck}.

\section{Proximity with respect to all extensions of $v$}
\label{S:one}

We work under the assumption that there exists a coherent good set of valuations $U_K$ on $K$. We first construct the set of local heights associated to the places in $U_K$ and then we define the global height. All our valuations in this section are normalized so that their value group is $\mathbb{Z}$.

For each finite extension field $L$ of $K$ and for each place $w$ of $L$ lying above a place in $U_K$, we let $\tilde{w}:L\rightarrow\mathbb{Z}_{\le 0}$ be defined as follows
$$\tilde{w}:=\min\{w,0\}.$$
Then the local height at $w$ of any element $x\in L$ is $\h_w(x):=-d(w)\tilde{w}(x)$. We define the global height of $x$ as 
$$\h(x):=\sum_{w}\h_w(x).$$
The above sum is a finite sum because there are finitely many $w$ such that $w(x)<0$ (see condition $(i)$ of Definition~\ref{D:good sets}).
Because $U_K$ is a coherent good set of valuations, the definition of the global height of an element $x$ does not depend on the particular choice of the field $L$ containing $x$ (see for example Chapter $4$ of \cite{teza}). The following two standard properties of the height will be used in our proof.

\begin{proposition}
\label{P:height inequalities}
For each $x,y\in K^{\alg}$, the following are true:

(i) $\h(xy)\le\h(x)+\h(y)$.

(ii) $\h(x+y)\le\h(x)+\h(y)$.
\end{proposition}

\begin{proof}
The proof is immediate using the definition of height and the triangle inequality for each valuation.
\end{proof}

For a point $P:=(x_1,\dots,x_g)\in\mathbb{G}_a^g(L)$, we define the local height of $P$ at a place $w$ of $L$ lying above a place in $U_K$, as follows:
$$\h_w(P):=\max\{\h_w(x_1),\dots,\h_w(x_g)\}.$$
Then the global height of $P$ is $\h(P):=\sum_w \h_w(P)$.

Next we define the heights associated to a Drinfeld module $\phi:A\rightarrow K\{\tau\}$ (see \cite{teza} for more details). We fix a non-constant $a\in A$. For each finite extension $L$ as above and for each place $w$ of $L$ as above, we define
$$V_w(x):=\lim_{n\rightarrow\infty}\frac{\tilde{w}(\phi_{a^n}(x))}{\deg(\phi_{a^n})},$$
for each $x\in L$.

Then the canonical local height of $x$ at $w$ with respect to $\phi$ is $\hhat_w(x):=-d(w)V_w(x)$. Finally, the canonical global height of $x$ with respect to $\phi$ is $\hhat(x):=\sum_w\hhat_w(x)$. By the same reasoning as in \cite{Den} (see part $3)$ of Th\'{e}or\`{e}me $1$) or in \cite{poo} (see part $(2)$ of Proposition $1$) we can show that there exists a positive constant $C_0$ such that for every $x\in K^{\alg}$, 
\begin{equation}
\label{E:difference heights}
|\h(x)-\hhat(x)|\le C_0.
\end{equation}
Moreover, the constant $C_0$ is easily computable in terms of $\phi$ (see \cite{poo}). 

For each point $P:=(x_1,\dots,x_g)\in\mathbb{G}_a^g(L)$ and for each place $w$ of $L$ as above, we define the canonical local height of $P$ at $w$ as $\hhat_w(P):=\max\{\hhat_w(x_1),\dots,\hhat_w(x_g)\}$. The canonical global height of $P$ is $\hhat(P):=\sum_w\hhat_w(P)$.

Using \eqref{E:difference heights} and Proposition~\ref{P:height inequalities} we prove the following result.
\begin{lemma}
\label{L:small height}
Let $L$ be a finite extension of $K$ and let $f\in L[X_1,\dots,X_g]$. There exists a constant $C(f)>0$ such that for every $P\in\mathbb{G}_a^g(K^{\alg})$, if $P$ is a torsion point for $\phi$, then $\h(f(P))\le C(f)$.
\end{lemma}

\begin{proof}
Using Proposition~\ref{P:height inequalities} $(i)$, it suffices to prove Lemma~\ref{L:small height} under the assumption that $f$ is a monomial. Hence, assume $f:=cX_1^{\alpha_1}\cdot\dots\cdot X_g^{\alpha_g}$ for some $c\in L$ and $\alpha_1,\dots,\alpha_g\in\mathbb{Z}_{\ge 0}$. Let $P=(x_1,\dots,x_g)$. We know that for each $i$, $x_i\in\phi_{tor}$. Hence $\hhat(x_i)=0$ for each $i$. Using \eqref{E:difference heights} we conclude that $\h(x_i)\le C_0$ for each $i$. Therefore, an application of Proposition~\ref{P:height inequalities} $(ii)$ concludes the proof of our Lemma~\ref{L:small height}.
\end{proof}

We proceed to the proof of Theorem~\ref{T:tv}.
\begin{proof}[Proof of Theorem~\ref{T:tv}.]
Let $f_1,\dots,f_m$ be a set of polynomials in $K[X_1,\dots,X_g]$ with integral coefficients at $v$, which generate the vanishing ideal of $X$. It suffices to prove that for each such polynomial $f_i$ and for every finite extension $L$ of $K$ and for every torsion point $P\in\mathbb{G}_a^g(L)$, either $f_i(P)=0$ or there exists a place $w|v$ of $L$ such that $w(f_i(P))\le\frac{C(f_i)}{d(v)}e(w|v)$, where $C(f_i)$ is the constant corresponding to $f_i$ as in Lemma~\ref{L:small height}. Then we obtain Theorem~\ref{T:tv} with $C:=\max_i\frac{C(f_i)}{d(v)}$.

Assume for some $i\in\{1,\dots,m\}$ and for some torsion point $P\in\mathbb{G}_a^g(L)$, $w(f_i(P))>\frac{C(f_i)}{d(v)}e(w|v)$ for every place $w|v$ of $L$. Then
\begin{equation}
\label{E:many zeros}
\sum_{w|v}d(w)\cdot w(f_i(P))>\frac{C(f_i)}{d(v)}\sum_{w|v}d(w)e(w|v)=\frac{C(f_i)}{d(v)}\sum_{w|v}\frac{d(v)f(w|v)e(w|v)}{[L:K]}=C(f_i)>0
\end{equation}
because $\sum_{w|v}f(w|v)e(w|v)=[L:K]$, as $v$ is a coherent valuation. If $f_i(P)\ne 0$, then \eqref{E:many zeros} yields that the set $S$ of places of $L$ lying above places in $U_K$ for which $f_i(P)$ is non-integral, is non-empty. Moreover, using \eqref{E:many zeros} and the sum formula for the nonzero element $f_i(P)\in L$, we conclude
\begin{equation}
\label{E:many poles}
\sum_{w\in S}d(w)\cdot w(f_i(P))<-C(f_i).
\end{equation}
Therefore, by the definition of the local heights we get
\begin{equation}
\label{E:big local heights}
\sum_{w\in S}\h_w(f_i(P))>C(f_i).
\end{equation}
Using the definition of the global height and \eqref{E:big local heights} we conclude $\h(f_i(P))>C(f_i)$. This last inequality contradicts Lemma~\ref{L:small height} because $P$ is a torsion point. This shows that $f_i(P)=0$ assuming $f_i(P)$ is close $w$-adically to $0$ for each $w|v$. This concludes the proof of our Theorem~\ref{T:tv}.
\end{proof}

\section{Proximity with respect to one fixed extension of $v$}
\label{S:future}

In this Section~\ref{S:future} we work under the hypothesis that the valuation $v$ of $K$ does not extend the valuation $v_{\infty}$ of $\Frac(A)$ in case $\phi:A\rightarrow K\{\tau\}$ is a Drinfeld module of generic characteristic. We also work with a fixed completion $K_v$ of $K$ at $v$ and with its algebraic closure $K_v^{\alg}$. In this section, the value group of our valuation $v$ is $\mathbb{Q}$, while its restriction to $K$ has value group $\mathbb{Z}$.

We first reduce Theorem~\ref{T:Mattuck} to the following Lemma~\ref{L:ball}.
\begin{lemma}
\label{L:ball}
Let $\phi:A\rightarrow K\{\tau\}$ be a Drinfeld module and let $v$ be a discrete valuation on $K$. If $\phi$ has generic characteristic, assume moreover that $v$ does not lie over the place $v_{\infty}$ of $\Frac(A)$. There exists a positive constant $C_v$ depending only on $\phi$ and $v$ such that in the ball $$\{x\in K_v^{\alg}\mid v(x)\ge C_v\}$$ there are no nonzero torsion points of $\phi$.
\end{lemma}
Lemma~\ref{L:ball} shows that for each place $v$ which does not lie over $v_{\infty}$ (if $\phi$ has generic characteristic), $\phi_{\tor}$ is discrete in the $v$-adic topology. If $\phi$ has finite characteristic, then $\phi_{\tor}$ is discrete with respect to each valuation $v$ (without any restriction). Moreover, as it will be shown in the proof of Lemma~\ref{L:ball}, the constant $C_v$ is easily computable in terms of $\phi$ and $v$.

\begin{proof}[Proof of Theorem~\ref{T:Mattuck}.]
We prove Theorem~\ref{T:Mattuck} using the result of Lemma~\ref{L:ball}. Let $Q:=(y_1,\dots,y_g)$. Let $\beta_i:=\max\{0,-v(y_i)\}$ for each $i\in\{1,\dots,g\}$. Let $\pi_v\in K$ be an uniformizer for $v$, i.e. $v(\pi_v)=1$. Then for each $i\in\{1,\dots,g\}$, the linear polynomial $\pi_v^{\beta_i}(X_i-y_i)\in L[X_1,\dots,X_g]$ has integral coefficients at $v$ and vanishes at $Q$.

We know (see Lemma $5.2.5$ of \cite{teza}) that there exists an absolute constant $M_v\le 0$ depending only on $\phi$ and $v$ such that for every torsion point $x\in\phi_{\tor}$, $v(x)\ge M_v$ (because otherwise, $x$ has positive local height at $v$, contradicting the fact that each local height of a torsion point is $0$). Then for each point $P:=(x_1,\dots,x_g)\in\phi_{\tor}^g$, if for some $i\in\{1,\dots,g\}$, $v(y_i)=-\beta_i<M_v\le v(x_i)$, then $v(x_i-y_i)=v(y_i)$. In this case, $\lambda_v(P,Q)\le v\left(\pi_v^{\beta_i}(x_i-y_i)\right)=0$. Therefore, in case for some $i\in\{1,\dots,g\}$, $v(y_i)<M_v$, we obtained an absolute upper bound for the $v$-adic distance of a torsion point to $Q$.

Assume from now on in this proof that for every $i\in\{1,\dots,g\}$, $v(y_i)\ge M_v$. Hence $\beta_i\le -M_v$. We compute the $v$-adic distance between a torsion point $P:=(x_1,\dots,x_g)\in\phi_{\tor}^g$ and $Q$. We obtain: $$\lambda_v(P,Q)\le\min_{i=1}^gv(\pi_v^{\beta_i}(x_i-y_i))=\min_{i=1}^g\left(\beta_i+v(y_i-x_i)\right)\le -M_v+\min_{i=1}^gv(x_i-y_i).$$
Therefore, in order to prove Theorem~\ref{T:Mattuck} it suffices to show that $$\min_{i=1}^gv(x_i-y_i)$$ is uniformly bounded from above when $(x_1,\dots,x_g)\in\phi_{\tor}^g\setminus\{(y_1,\dots,y_g)\}$. But Lemma~\ref{L:ball} shows that for each $i$, there is at most one torsion point of $\phi$ in the ball
\begin{equation}
\label{E:ball_displacement}
\{x\in K_v^{\alg}\mid v(x-y_i)\ge C_v\},
\end{equation}
because otherwise there would be at least one nonzero torsion point of $\phi$ in $\{x\in K_v^{\alg}\mid v(x)\ge C_v\}$ after translating the ball in \eqref{E:ball_displacement} by a torsion point of $\phi$ which lies inside the ball from \eqref{E:ball_displacement}. Therefore, $\lambda_v(P,Q)$ is indeed uniformly bounded from above for $P\in\phi_{\tor}^g\setminus\{Q\}$ because there is at most one torsion point $P\in\phi_{\tor}^g$ such that $\lambda_v(P,Q)>-M_v+C_v$.
\end{proof}

We proceed to the proof of Lemma~\ref{L:ball}.
\begin{proof}[Proof of Lemma~\ref{L:ball}.]
We first choose $t\in A$ satisfying certain properties according to the two cases we have: $\phi$ has generic characteristic or not.

\emph{Case (i)}. $\phi$ has generic characteristic.

Let $\pp$ be the nonzero prime ideal of $A$ which is contained in the maximal ideal of the valuation ring of $v$ (we are using the fact that $v$ does not lie over $v_{\infty}$ to derive that all the elements of $A$ are integral at $v$). We fix $t\in\pp\setminus\{0\}$.

\emph{Case (ii)}. $\phi$ has finite characteristic.

Let $\pp$ be the characteristic ideal of $\phi$. By the hypothesis for our \emph{Case (ii)}, $\pp$ is nonzero. We fix $t\in\pp\setminus\{0\}$.

Let $\phi_t=\sum_{i=r_0}^ra_i\tau^i$, where $a_{r_0}\ne 0$. In finite characteristic, $r_0\ge 1$, while in generic characteristic, $r_0=0$ and $v(a_0)\ge 1$ (by our choice of $t$). We let $C_v$ be the smallest positive integer larger than all of the numbers from the following set:
$$S:=\{-\frac{v(a_{r_0})}{q^{r_0}-1}\}\cup\{\frac{v(a_{r_0})-v(a_i)}{q^i-q^{r_0}}\mid r_0<i\le r\}.$$
We note that if $\phi$ has generic characteristic, then $r_0=0$ and so, $q^{r_0}=1$. Then the denominator of the first fraction contained in $S$ is $0$. So, because the numerator $-v(a_0)\le -1$, that fraction equals $-\infty$ and so, any integer is larger than it, i.e. if $\phi$ has generic characteristic, we may disregard the first fraction in the definition of $S$. As we will see in our proof, that first fraction will only be used in the finite characteristic case.

\begin{Claim}
\label{C:decisive}
If $x\in K_v^{\alg}\setminus\{0\}$ satisfies $v(x)\ge C_v$, then $v(\phi_t(x))=v\left(a_{r_0}x^{q^{r_0}}\right)>v(x)\ge C_v$. In particular, $\phi_t(x)\ne 0$.
\end{Claim}

\begin{proof}[Proof of Claim~\ref{C:decisive}.]
Because $v(x)\ge C_v$, then for each $i\in\{r_0+1,\dots,r\}$
\begin{equation}
\label{E:t-val}
v\left(a_ix^{q^i}\right)>v\left(a_{r_0}x^{q^{r_0}}\right).
\end{equation}
Inequality \eqref{E:t-val} shows that $v\left(\phi_t(x)\right)=v\left(a_{r_0}x^{q^{r_0}}\right)$. In particular, this shows $\phi_t(x)$ does not equal $0$, because its valuation is not $+\infty$. Hence 
\begin{equation}
\label{E:t-ineq}
v(\phi_t(x))=v(a_{r_0})+q^{r_0}v(x).
\end{equation}
If $\phi$ has generic characteristic, then \eqref{E:t-ineq} shows that $v(\phi_t(x))=v(a_0)+v(x)\ge v(x)+1>C_v$. If $\phi$ has finite characteristic, then using that 
$$v(x)\ge C_v>-\frac{v(a_{r_0})}{q^{r_0}-1}$$
we conclude $v(\phi_t(x))=v(a_{r_0})+q^{r_0}v(x)>v(x)\ge C_v$. 
\end{proof}

Claim~\ref{C:decisive} shows that for every nonzero $x\in K_v^{\alg}$ satisfying $v(x)\ge C_v$, the sequence $\{v(\phi_{t^n}(x))\}_{n\ge 0}$ is strictly increasing. Hence, $x\notin\phi_{\tor}$, because if $x$ were torsion, then the sequence $\{\phi_{t^n}(x)\}_{n\ge 0}$ would contain only finitely many distinct elements. This concludes the proof of Lemma~\ref{L:ball}.
\end{proof}

\address{Dragos Ghioca, Department of Mathematics, McMaster University, Hamilton Hall, Room 218,
Hamilton, Ontario L8S 4K1, Canada}

\email{dghioca@math.mcmaster.ca}

\end{document}